\documentclass[12pt]{article}
\newcommand{\lb}[1]{\label{#1}}\newcommand{\reff}{\ref}\newcommand{\rf}[1]{(\reff{#1})}\newcommand{\cit}[1]{\cite{#1}}\newcommand{\bibi}[1]{\bibitem{#1}}
\newcommand{\longue}[1]{}
\usepackage{amssymb}
\usepackage{amsmath}
\usepackage[latin1]{inputenc}

\newtheorem{lm}{{\bf Lemma}}
\newtheorem{theor}[lm]{{\bf Theorem}}
\newtheorem{deff}[lm]{{\bf Definition}}

\newtheorem{cj}[lm]{{\bf Conjecture}}

\newtheorem{cor}[lm]{{\bf Corollary}}
\newtheorem{exemple}[lm]{{\bf Example}}
\newtheorem{question}{{\bf Question}}

\newtheorem{prop}[lm]{{\bf Proposition}}
\newcommand{\sep}{$\!${\rm.}\ }

\newcommand{\qu}[2]{\begin{question}{\lb{#1}}\sep{\sl #2}\end{question}}

\newcommand{\theo}[2]{\begin{theor}{\lb{#1}}\sep{\sl 
#2}\end{theor}}

\newcommand{\coro}[2]{\begin{cor}{\lb{#1}}\sep{\sl #2}\end{cor}}
\newcommand{\propo}[2]{\begin{prop}{\lb{#1}}\sep{\sl #2}\end{prop}}

\newcommand{\eq}[2]{\begin{equation}#2\ \ \  \ \lb{#1}\end{equation}}

\newcommand{\pr}{\noindent {\sl Proof\;\sep }}
\newcommand{\ep}{\hfill \framebox[2mm]{\ } \medskip}
\newcommand{\tq}{\ ;\ }
\newcommand{\be}{\begin{enumerate}}
\newcommand{\ee}{\end{enumerate}}
\newcommand{\noi}{\noindent}
\newcommand{\med}{\medskip}

\newcommand{\ts}{\textstyle}
\newcommand{\ds}{\displaystyle}
\newcommand{\lp}{\left(}
\newcommand{\rp}{\right)}

\renewcommand{\bar}{\overline}

\newcommand{\ssi}{\Leftrightarrow}

\newcommand{\F}{{\mathcal F}}

\newcommand{\f}{\varphi}

\renewcommand{\S}{{\mathcal S}}

\newcommand{\N}{\mathbb{N}}

\newcommand{\ind}{{\rm Ind}\,}

\newcommand{\np}{\newpage}
\newcommand{\ed}{\end{document}}

\newcommand{\wlg}{Without loss of generality}
\title{
{\Large 
What is the minimal cardinal of a family\\ 
which shatters all $d$-subsets of a finite set?
}
}
\author{
N.\ Chevallier and A.\ Fruchard}
\date{
\today}
\begin{document}
\maketitle
%
%
\noi In this note, $d\leq n$ are positive integers.
Let $S$ be a finite set of cardinal $|S|=n$ and let $2^S$ denote 
its power set, i.e. the set of its subsets.
A {\sl $d$-subset of $S$} is a subset of $S$ of cardinal $d$.
Let $\F\subseteq2^S$ and $A\subseteq S$.
The {\sl trace of $\F$ on $A$} is the family $\F_A=\{E\cap A\tq E\in\F\}$. 
One says that $\F$ {\sl shatters} $A$ if $\F_A=2^A$. 
The {\sl VC-dimension of $\F$} is the maximal cardinal of a subset of $S$
that is shattered by $\F$~\cit{vc}.
The following is well-known~\cit{vc,s,sh}: 
\theo{t1}{
{\rm(Vapnik-Chervonenkis, Sauer, Shelah)}\\
If VC-dim$(\F)\leq d$ (i.e. if $\F$ shatters no $(d+1)$-subset of $S$)
then $|\F|\leq c(d,n)$, where 
$$
c(d,n)=\ts\binom n0+\dots+\binom nd.
$$
Moreover this bound is tight: 
It is achieved e.g. for $\F=\binom S{\leq d}$, the family of
all $k$-subsets of $S$, $0\leq k\leq d$.
}
A first natural question is:
\qu{q1}{
Assume a family $\F\subseteq2^S$ is maximal for the inclusion among all
families of VC-dimension at most $d$. Does $\F$ always have the maximal
possible cardinal $c(d,n)$?
}
Let us define the {\sl index of} $\F$ as follows:
$$
\ind\F=\max\{d\in\{0,...,n\}\tq\F\mbox{ shatters all {\it d}-subsets of }S\}.
$$
Let $C(d,n)=\min\{|\F|\tq\ind\F=d\}$.
For instance, we have $C(1,n)=2$, with the (only possible) choice 
$\F=\{\emptyset,S\}$.
Of course we have $2^d\leq C(d,n)\leq2^n$.
The question is:
\qu{q2}{
Give the exact value of $C(d,n)$ for $2\leq d\leq n$. 
If this is not possible, give lower and upper bounds as accurate as possible.
}
A well-known duality yields another formulation of Question~\reff{q2}.
Let $\f:S\to2^\F,\;a\mapsto\{E\in\F\tq a\in E\}$ and set $\S=\f(S)$.
In this manner, we have for all $a\in S$ and all $E\in\F$:
\eq2{
a\in E\ssi E\in\f(a).
}
One can check that $\F$ shatters $A\subseteq S$ if and only if,
for every partition $(B,C)$ of $A$ (i.e. $A=B\cup C$ and $B\cap C=\emptyset$)
the intersection 
$\Big(\bigcap_{b\in B}\f(b)\Big)\cap\Big(\bigcap_{c\in C}\bar{\f(c)}\Big)$
is nonempty, where the notation $\bar Y$ stands for $\F\setminus Y$.

If $\ind\F\geq2$, then $\f$ is a one-to-one correspondance from $S$ 
to $\S$, hence we have $\log n\leq{C(d,n)}$ for all $2\leq d\leq n$, where
$\log$ denotes the logarithm in base $2$.
\bigskip
\np
\noi{\large\bf The case $d=2$}.
Using for instance the binary expansion, it is easy to show that 
the order of magnitude of $C(2,n)$ is actually $\log n$.
The next statement refines this.
%
\propo{p2}{
If $n=\frac12\binom{2l}l=\binom{2l-1}{l-1}$, then $C(2,n)=2l$.
}
\pr
(Recall the notation $\bar A=\F\setminus A$.)
We first prove by contradiction that $C(2,n)>2l-1$.
Actually, if a family $\F$ of 
subsets of $S$ 
shatters all 2-subsets of $S$, then 
the image $\S\subseteq2^\F$ of $S$ by $\f$ must satisfy
\eq1{
\forall A\neq B\in\S,\;A\cap B,\;A\cap\bar B,\;\bar A\cap B,
\mbox{ and }\bar A\cap\bar B\mbox{ are nonempty}.
}
In particular $\S$ is a Sperner family of $\F$
(i.e. an antichain for the partial order of inclusion; one finds 
several other expressions in the literature: `Sperner system', 
`independent system', `clutter', `completely separating system', etc.).
For a survey on Sperner families and several generalizations, we refer e.g.
to~\cit b and the references therein.

Assume now that $|\F|=2l-1$; it is known~\cit{sp,k,l} that 
all Sperner families of $\F$ have a cardinal at most $\binom{2l-1}{l-1}$,
and that there are only  two Sperner families of maximal cardinal:
the families $\binom\F{l-1}$ and $\binom\F l$, i.e. of $(l-1)$-subsets, 
resp. $l$-subsets  of $\F$.
\longue{
This can be seen using the so-called Lubell-Yamamoto-Meshalkin
(LYM) inequality. The LYM inequality states that,
if $a_k$ denotes the number of sets of size $k$ in a Sperner family of a $n$-set,
then $\ds\sum_{k=0}^n\frac{a_k}{\binom nk}\leq1$.
}
However, none of these families satisfies both $A\cap B$ and $\bar A\cap\bar B$
nonempty in~\rf1.
As a consequence, we must have $|\F|\geq2l$.
 
Conversely, let $S=\{a_1,\dots,a_n\}$,
consider $\binom{\{1,\dots,2l\}}l$,
the set of $l$-subsets of $\{1,\dots,2l\}$, and choose one element in each pair of 
complementary $l$-subsets.
We then obtain a family $\{A_1,\dots,A_n\}$ which satisfies~\rf1.
Now we set $\F=\{E_1,\dots,E_{2l}\}$, with $E_i=\{a_j\tq i\in A_j\}$.
The characterization~\rf2 shows that $\F$ shatters every $2$-subset of $S$.
\ep

The proof of the following statement is straightforward.
\coro{c3}{
If $\binom{2l-1}{l-1}<n\leq\binom{2l+1}l$, then $2l\leq C(2,n)\leq2l+2$.
}
The upper bound can be slightly improved:
One can prove that, if $\binom{2l-1}{l-1}<n\leq\binom{2l}{l-1}$, 
then $2l\leq C(2,n)\leq2l+1$.
\longue{
Actually, if $n=\binom{2l}{l-1}$, with the notation $k=2l+1$ and 
$\F=\{1,\dots,k\}$, the family 
$\S=\Big\{E\cup\{k\}\tq E\in\binom{\F\setminus\{k\}}l\Big\}$
shows that $C(2,n)\leq2l+1$.
}
\qu{q3}{
It seems that we have $C(2,n)=k$ if and only if 
$\binom{k-2}{\lfloor(k-1)/2\rfloor-1}<n\leq\binom{k-1}{\lfloor k/2\rfloor-1}$, where
$\lfloor x\rfloor$ denotes the integer part of $x$.
Is it true? Is it already known?
}
The first values are $C(2,2)=C(2,3)=4$, $C(2,4)=5$, $C(2,5)=\dots=C(2,10)=6$.
\longue{
{\noi\sl Proof of $C(2,5)=\dots=C(2,10)=6$.}
\med

Assume by contradiction that $C(2,5)=5$. Denote $S=\{a,b,c,d,e\}$, 
$\F=\{1,2,3,4,5\}$, ans $\S=\{A,B,C,D,E\}$.
Since we must have  $X\cap Y$, $X\cap\bar Y,$ $\bar X\cap Y$,  
and $\bar X\cap\bar Y$ nonempty for all $X,Y\in\S$,
we have 
\eq6{
\forall X\in\S\quad2\leq X\leq 3.
}
\wlg, we assume that at least three elements of $\S$, say $A,B,C$ are of 
cardinal $2$. If $A\cap B\cap C=\emptyset$, since 
\eq7{
A\cap D,\;B\cap D,~\mbox{ and }~C\cap D~\mbox{  are nonempty},
}
this forces $D$ to contain one of $A,B,$ or $C$, a contradiction, therefore
$A\cap B\cap C\neq\emptyset$.
\wlg, assume  $A=\{1,2\}$, $B=\{1,3\}$, $C=\{1,4\}$.
If $1\notin D\cup E$, then \rf7 yields $\{2,3,4\}\subseteq D$ and similarly
$\{2,3,4\}\subseteq E$, yielding $D=E$ by \rf6, a contradiction.
\wlg, assume $1\in D$, then since $A,B,C\not\subseteq D$,
this forces $D=\{1,5\}$, hence $1\notin E$ (otherwise $E$ contains one
of $A,B,C$, or $D$), hence $E=\{2,3,4,5\}$, 
hence $\bar A\cap\bar E\neq\emptyset$, a contradiction.
This gives $C(2,5)>5$. Since $C(2,n)$ is nondecreasing in $n$ and $C(2,10)=6$
by Proposition~\reff{p2}, the proof is complete.
\ep
}
Computer seems to be useless, at least for a naive treatment.
Already in order to obtain $C(2,11)=7$, 
we would have to verify that $C(2,11)>6$, i.e. to find, for each of the
$\binom{2^{11}}6\approx10^{17}$ 
families  $\F$ in $2^S$ some $2$-subset that is not shattered by the family.
(Alternatively, in the dual statement, we have to check ``only'' 
$\binom{2^6}{11}\approx7.10^{11}$ families $\S$ in $2^\F$.)
\bigskip

\noi{\large\bf The case $d\geq3$}.
From now, we assume $n\geq4$. 
\propo{p4}{
For all $3\leq d<n$, we have 
$C(d,n)\leq\frac{2^d}{d!}\,(3\log n)^d$.
}
The constant $3$ can be improved. The proof below shows that, for all
$a>1$ and all $n$ large enough, 
$C(d,n)\leq\frac{2^d}{d!}\,(a\log n)^d$.
\med

\pr
Let $\F_0\subset2^S$ be a minimal 
separating system of $S$, i.e. such that, for all $a,b\in S$
there exists $E_a^b\in\F_0$ which satisfies $b\notin E_a^b\ni a$. 
Since this amounts to choosing $\F_0$ minimal such that $\S=\f(S)$ is a
Sperner family for $\F_0$, we know that $|\F_0|=N$ if and only if 
$\binom{N-1}{\lfloor(N-1)/2\rfloor}<n\leq\binom{N}{\lfloor N/2\rfloor}$,
hence $N:=|\F_0|
\leq2+\log n+\frac12\log\log n
\leq3\log n$ 
since $n\geq4$. We assume $N\geq2$ in the sequel.
Given two disjoint subsets $B$ and $C$ of $S$ such that $|B\cup C|=d$, the set 
$
E_B^C=
\bigcap_{c\in C}\lp\bigcup_{b\in B}E_b^c\rp
$
contains $B$ and does not meet $C$.
Let $\F$ be the collection of all such sets 
$E_B^C$; then $\F$ shatters all subsets of $S$ of cardinal at most $d$.

To estimate $|\F|$, we consider $\F_k$ the collection of all such sets 
$E_B^C$, with $|B|=k$ (and thus $|C|=d-k$).
We have $|\F_k|=\binom Nk\binom{N-k}{d-k}$ (with $N=|\F_0|$).
Then we choose $\F=\bigcup_{k=0}^d\F_k$.
We obtain $|\F|\leq\sum_{k=0}^d\binom Nk\binom{N-k}{d-k}=\binom Nd2^d
\leq\frac{2^d}{d!}N^d\leq\frac{2^d}{d!}\,(3\log n)^d$.
%

\ep
\qu{q4}{
Is $(\log n)^{\lfloor d/2\rfloor\lfloor(d+1)/2\rfloor}$ the right order of
magnitude for $C(d,n)$?
}
By constructing auxiliary Sperner families from $\S$, it is possible to give
a better lower bound for $C(d,n)$ than only $C(d,n)\geq C(2,n)$.
For instance, in the case $d=3$, for all distinct $A,B,C\in\S$,
we must have $A\cap B\not\subseteq C$. One can check that this implies that
the family $\{A\cap B\tq A,B\in\S\}$ 
is a Sperner family, therefore we obtain 
$\binom n2\leq\binom{C(3,n)}{\lfloor C(3,n)\rfloor/2}$. 
Unfortunately, this does not modify the order of magnitude.
Already in this case $d=3$, we do not know whether $C(3,n)$ is of order $\log n$, 
$(\log n)^2$, or an intermediate order of magnitude.
Another formulation is:
\qu{q5}{
Prove or disprove: There exists $C>0$ such that, for all $k\in\N$, 
if $\F$ is a finite set of cardinal $k$ and $\S\subseteq2^\F$ satisfies
$\forall A,B,C\in\S,\;A\cap B\not\subseteq C$, 
then $|\S|\leq C\,2^{C\sqrt k}$.
}

\end{document}